\documentclass{article}
\usepackage{amsmath,amssymb,bm,bbm,mathrsfs,amscd}
\usepackage{color}
\usepackage{amsthm}
\usepackage{dsfont}
\usepackage{graphicx}
\usepackage{epsfig}
\usepackage{subfigure}

\def\qed{\hfill \vrule height 7pt width 7pt depth 0pt\medskip}
\def\beq{\begin{equation}}
\def\eeq{\end{equation}}

\newtheorem{theorem}{Theorem}

\newtheorem{proposition}[theorem]{Proposition}
\newtheorem{lemma}[theorem]{Lemma}
\newtheorem{corollary}[theorem]{Corollary}
\newtheorem{example}{Example}
\theoremstyle{remark}

\newcommand{\ba}{\begin{array}}
\newcommand{\ea}{\end{array}}

\newcommand{\be}{\begin{equation}}
\newcommand{\ee}{\end{equation}}

\newcommand{\mc}{\mathcal}

\newcommand{\1}{\mathbbm{1}}

\newcommand{\R}{\mathbb{R}}
\newcommand{\N}{\mathbb{N}}

\renewcommand{\P}{\mathbb{P}}

\DeclareMathOperator*{\argmax}{argmax}

\def\1{\mathds{1}}

\def\N{\mathbb{N}}

\def\R{\mathbb{R}}

\def\P{\mathbb{P}}

\newcommand{\RR}{{\mathbb{R}}}

\title{\LARGE \bf A game theoretic approach to a network allocation problem}

\author{Fabio Fagnani$^{1}$,  Barbara Franci$^{2}$
\thanks{$^{1}$Fabio Fagnani is with the Department of Mathematical Sciences, Politecnico di Torino, 
        {\tt\small fabio.fagnani@polito.it}}%
\thanks{$^{2}$Barbara Franci is with the Department of Mathematical Sciences, Politecnico di Torino,
        {\tt\small barbara.franci@polito.it}}%
}

\begin{document}

\maketitle
\thispagestyle{empty}
\pagestyle{empty}

\begin{abstract}
In this paper, we consider a network allocation problem  
motivated by peer-to-peer cloud storage models. The setting is that of a network of units (e.g. computers) that collaborate and offer each other space for the back up of the data of each unit.
We formulate the problem as an optimization problem, we cast it into a game theoretic setting and we then propose a decentralized allocation algorithm based on the log-linear learning rule. Our main technical result is to prove the convergence of the algorithm to the optimal allocation. We also present some simulations that show the feasibility of our solution and corroborate the theoretical results.
\end{abstract}


\section{Introduction}\label{sec:introduction}

Recently, cooperative storage cloud models based on peer-to-peer architectures have been proposed as valid alternatives to traditional centralized cloud storage services. The idea is quite simple: instead of using dedicated servers for the storage of data, the participants themselves offer space available on their connected devices to host data from other users. In this way, each participant has two distinct roles: that of a unit that needs external storage to securely back up its data, and that of a resource available for the back up of data of other users. This approach has in principle a number of relevant advantages with respect to traditional cloud storage models. First, it eliminates the need for a significant dedicated hardware investment so that the service should be available at (order of magnitude) lower cost. Second, it overcomes the typical problems related to the use of a single external provider as security threats, including man-in-the-middle attacks, and malware, or fragility with respect to technical failures.

On the wake of the successful peer-to-peer file sharing model of applications like BitTorrent and its lookalike, the same philosophy may well be leveraged on a different but very similar application service like storage. Indeed a slew of fledgling and somehow successful startups are entering in this market niche. Among the most noteworthy examples are the platforms Sia \text{http://sia.tech} and Storj \text{https://storj.io} (see for some information \cite{SIA} and \cite{STORY}).

Clearly, a completely decentralized peer-to-peer model must account for some challenging technical difficulties that are absent in a centralized cloud model. Firstly, security and privacy must be carefully implemented by ensuring end-to-end encryption resistant to attackers, as well suitable coding to insure recovering from failure of some units. In addition, the model must account for the latency, performance, and downtime of average user devices. Albeit the above technical issues are challenging, they can be addressed with the right tools and architectures available at current state of the art technology and will not be considered in this paper. 

A core part of such cooperative storage model is the mechanism by which users are made to interact, collaborate and share their storage commodity with each other. In the existing platforms as Sia and Storj, this part is worked out at the level of a central server to which all units are connected. At the best of our knowledge, the design of such a mechanism in a decentralized fashion 
 has not been yet theoretically addressed and studied in the literature. Our contribution goes in this direction.

In this paper we present, accompanied by a rigorous mathematical analysis, a cooperative fully distributed algorithm through which a network of units (e.g. computers) can collaborate and offer each other space for the back up of the data of each unit. In this model, there is no need for central supervision and it can easily incorporate features that we want the system to possess depending on the application, as for instance, enforcing structure on the way data of each unit is treated (aggregate or rather disgregate in the back up process), avoiding congestion phenomena in the use of the resources, differentiate among resources on the basis of their reliability. A different version of the algorithm and without analytical results was presented in \cite{GTA}.

We formulate the problem as an optimal network allocation problem where a population of units is connected through a graph and each of them possesses a number of items that need to be allocated among the neighboring units. Each unit in turn offer a certain amount of storage space where neighboring units can allocate their items. The optimal allocation is the one maximizing a given functional that depends on the allocation status of each unit and that incorporates the desired features we want the solution to possess. 

In order to solve the optimization problem in a scalable decentralized fashion we cast the allocation problem into a game theoretic framework and we design the algorithm using a learning dynamics. 
The use of game theory to solve distributed optimization problems and, more in general, in the design and control of large scale networked engineering systems is becoming increasingly popular \cite{GTW, GTL, GTDS, DGDO}. The basic idea is that of modeling system units as rational agents whose behavior consists in a selfish search for their maximum utility. The goal is to design both the agents utility functions and a learning adaptation rule in such a way that the resulting global behavior is (close to) the desired one, the maximum in the specific case of optimization problems. Two are the challenges that typically we need to face. First, design agent utility functions that only uses information present at the level of the single units and that leads to a game whose Nash equilibria contain the desired configurations. Second, design the learning mechanism through which system converges to a desired Nash. 

For optimization purposes, an interesting strategy \cite{DGDO, DWG} is to design utility functions so to yield a potential game whose potential coincide with the reward functional of the problem and then consider the log-linear learning dynamics (also known as noisy best response) \cite{SMS, GTDS, LOG}. Under certain assumptions, this rule is known to lead to a time-reversible ergodic Markov chain whose invariant probability distribution is a Gibbs measure with energy function described by the potential and (for small noise parameter) has its peak on the maxima of the potential. In this paper we follow this road. 

For a very general family of functionals having an additive separable form, namely that can be expressed as sums of terms depending on the various units, we define a game by setting the utility function of each unit as simply the sum of those addends in the functional involving the unit itself and its neighbors, while the action set of a unit consists of the vectors describing the allocation among its various neighbors. The game so defined is easily shown to be potential with potential given by the original functional. The game, however, possesses a key critical feature: because of the hard storage constraints of the various resources, units are not free to choose their actions as they want, but they are constrained from the choice made by other units. For instance, if a unit is saturating the space available in a certain resource, other units connected to the same resource will not be able to use it. In cooperative cloud storage models where resources are common users, this hard storage constraint is a very natural assumption and can not be relaxed. This property is non-classical in game theory and has remarkable consequences on the structure of Nash equilibria and the behavior of the best response dynamics that is not guaranteed in general to approximate the optimum. Indeed, constrained equilibrium problems and convergence algorithm are widely studied in literature and they are known as Generalized Nash Equilibrium Problems (\cite{GNEP,GenPot} and reference therein). 

The main technical contribution of this paper is to show that, despite these hard constraints, under mild technical conditions, a family of dynamics having their core on the log-linear learning rule converge to the desired solution. More precisely, by a careful analysis of the connectivity properties of the transition graph associated to the Markov process, we will obtain two results: (i) if there is enough space for allocation to take place, under the proposed algorithms, all units will complete their allocation in finite time with probability one; and (ii), under a slight stronger assumption, the allocation configuration converges, when the noise parameter approaches $0$, to a maximum of the original functional. At the best of our knowledge, this analysis is new in game theory.

We want to remark that, the type of functionals considered are typically non-convex (even when relaxed to continuous variables) so that many algorithms for distributed optimization may fail to converge to the global maximum. In addition, the proposed algorithm presents a number of interesting features. The algorithm is decentralized and adapted to any predefined graph. For the cloud application we have in mind, the choice of the graph topology can be seen as a design parameter that allows to control the computational complexity at the units level. It is asyncronous and it is robust with respect to temporary disconnection of units. Moreover, it is intrinsically open-ended: if new data or new units enters into the system, a new run of the algorithm will automatically permit the allocation of the new data and, possibly, the redistribution of the data stored by the old units to take advantage of new available space.

The remaining part of this section is dedicated to some literature review. In Section \ref{model} we formally define the network allocation problem and recall some basic facts proven in \cite{GTA} (in particular, a necessary and sufficient condition for the allocation problem to be solvable). We then introduce a family of functionals and define the optimal allocation problem. Section \ref{sec: algorithm} is devoted to cast the problem to a potential game theoretic framework \cite{ROS,Pot} and to propose a distributed algorithm that is an instance of a noisy best response dynamics. The main technical part of the paper is Section \ref{sec: analysis} where the fundamental results Theorem \ref{theo main 1} and Corollary \ref{cor main 2} are stated and proven. Theorem \ref{theo main 1} ensures that the algorithm  reaches a complete allocation with probability one, if a complete allocation is indeed possible. Corollary \ref{cor main 2} studies the asymptotic behavior of the algorithm and explicitly exhibits the invariant probability distribution. Consequence of Corollary \ref{cor main 2} is that in the double limit when time goes to infinity and the noise parameter goes to $0$, the algorithm converges to a Nash equilibrium that is, in particular, a global maximum of the potential function. This guarantees that the solution will indeed be close to the global welfare of the community. Finally, Section \ref{sec simulation} is devoted to the presentation of a set of simulations that show a practical implementation of the algorithm. Though we work out relatively simple examples, our simulations show the good properties of the algorithm, its scalability properties in terms of speed and complexity, and illustrate the effect of the parameters of the utility functions in the solution reached by the algorithm. A conclusion section ends the paper.

\subsection{Related Work}
Our problem fits into the wide class of distributed resource allocation problems. Among the many applications where such problems arise in a similar form to the one proposed in this paper, we can cite
cloud computing \cite{CRAG}, network routing \cite{SRPA}, vehicle target assignment \cite{VTA}, content distribution \cite{MSG}, graph coloring \cite{GC}. The game theoretic approach to allocation problems and the consequent design of distributed algorithms has been systematically addressed in \cite{DWG, DGDO, GTDS, LOG} where general techniques for the choice of the utility functions and of the dynamic learning rule have been proposed. 

The model proposed in this paper and the algorithm based on noisy best response dynamics, is inspired by this literature. A key aspect of our model and that makes it different from the models treated in the above literature is the fact that resources have hard storage limitations. This is a natural feature of the distributed cloud storage problem considered in this paper and that, to our knowledge, had not yet been previously analyzed.

\section{The cooperative storage model}\label{model}

Consider a set $\mc X$ of units that play the double role of users who have to allocate externally a back up of their data, as well resources where data from other units can be allocated. Generically, an element of $\mc X$ will be called a unit, while the terms user and resource will be used when the unit is considered in the two possible roles of, respectively, a source or a recipient of data. We assume units to be connected through a directed graph $\mc G=(\mc X,\mc E)$ where a link $(x,y)\in\mc E$ means that unit $x$ is allowed to store data in unit $y$. We denote by $$N_x:=\{y\in\mc X\,|\, (x,y)\in\mc E\},\quad N^-_y:=\{x\in\mc X\,|\, (x,y)\in\mc E\}$$
respectively, the out- and the in-neighborhood of a node. Note the important different interpretation in our context: $N_x$ represents the set of resources available to unit $x$ while $N^-_y$ is the set of units having access to resource $y$. If $D\subseteq \mc X$, we put $N(D)=\cup_{x\in D}N_x$ and $N^-(D)=\cup_{x\in D}N_x^-$. 

We imagine the data possessed by the units to be quantized atoms of the same size. Each unit $x$ is characterized by two non negative integers:
\begin{itemize}
\item $\alpha_x$ is the number of data atoms that unit $x$ needs to back up into his neighbors, 
\item $\beta_x$ is the number of data atoms that unit $x$ can accept and store from his neighbors.
\end{itemize}
The numbers $\{\alpha_x\}$ and $\{\beta_x\}$ will be assembled into two vectors denoted, respectively, $\alpha$ and $\beta$. 
Given the triple $(\mc G, \alpha, \beta)$, we define a \emph{partial state allocation}
as any matrix $W\in\N^{\mc X\times \mc X}$ that satisfies the following conditions
\begin{enumerate}
\item[(P1)]  $W_{xy}\geq 0$ for all $x,y$ and $W_{xy}=0$ if $(x,y)\not\in \mc E$.
\item[(P2)]  $W^x:=\sum\limits_{y\in\mc X}W_{xy}\leq\alpha_x$ for all $x\in\mc X$.
\item[(P3)] $W_y:=\sum\limits_{x\in\mc X}W_{xy}\leq \beta_y$ for all $y\in\mc X$.
\end{enumerate}
We interpret $W_{xy}$ as the number of pieces of data that $x$ has allocated in $y$ under $W$. Property (P1) enforces the graph constraint: $x$ can allocate in $y$ iff $(x,y)\in \mc E$. Property (P2) says that a unit can not allocate more data than the one it owns, and, finally, (P3) describes the storage constraint at the level of units considered as resources. Whenever $W$ satisfies (P2) with equality for all $x\in\mc X$, we say that $W$ is an \emph{allocation state}. The set of partial allocation states and the set of allocation states are denoted, respectively, with the symbols $\mc W_p$ and $\mc W$. 
We will say that the allocation problem is solvable if a state allocation $W\in\mc W$ exists. 

\subsection{Existence of allocations}
The following result gives a necessary and sufficient condition for the existence of allocations. The proof, which follows from Hall's theorem, can be found in \cite{GTA}. 

\begin{theorem}\label{theo allocation exists} Given $(\mc G, \alpha, \beta)$, there exists a state allocation iff the following condition is satisfied:
\beq\label{cond allocation exists}
\sum\limits_{x\in D}\alpha_x\leq \sum\limits_{y\in N(D)}\beta_y\quad\forall D\subseteq \mc X
\eeq
\end{theorem}

We first analyze the existence of allocations in the simple case of a complete network.


\begin{example}\label{ex: complete all} If $\mc G$ is complete, we have that $N(\{x\})=\mc X\setminus\{x\}$ while $N(D)=\mc X$ for all $D$ such that $|D|\geq 2$. 
Therefore, in this case,
condition (\ref{cond allocation exists})  reduces to
\beq\label{cond allocation exists complete}
\alpha_x\leq \sum\limits_{y\neq x}\beta_y,\; \forall x\in\mc X\qquad 
\sum\limits_{x\in\mc X}\alpha_x\leq \sum\limits_{y\in \mc X}\beta_y
\eeq
\end{example}

We now focus on the special but intersting case when all units have the same amount of data to be stored and the same space available, namely, $\alpha_x=a$, $\beta_x=b$ for every  $x\in\mc X$. 
In this case, condition (\ref{cond allocation exists complete}) that characterizes the existence of allocations for the complete graph, simply reduces to $a\leq b$.  In this case, among the possible allocation states there are those where each unit uses only one resource: given any permutation $\sigma :\mc X\to\mc X$ without fixed points, we can consider
\beq\label{extremal Nash}W^\sigma_{xy}=\left\{\begin{array}{ll} a &\;{\rm if} \, \sigma(x)=y\\
0 &\;{\rm otherwise} \end{array}\right.
\eeq
In general, an allocation state as $W^{\sigma}$ in (\ref{extremal Nash}) of the example above where each unit uses just one resource and each resource is only used by one unit, is called a \emph{matching allocation state}. Existence of matching allocation states is guaranteed for more general graphs than the complete ones.

\begin{proposition}\label{prop: existence matching all} Let $\mc G=(\mc X,\mc E)$ be any graph and assume that $\alpha_x=a$, $\beta_x=b$ for every  $x\in\mc X$. The following conditions are equivalent:
\begin{enumerate}
\item[(i)] There exists an allocation state $W$;
\item[(ii)] There exists a matching allocation state;
\item[(iii)] $a\leq b$ and $|D|\leq |N(D)|$ for every subset $D\subseteq \mc X$.
\end{enumerate}
\end{proposition}
\begin{proof}
(ii) $\Rightarrow$ (i) is trivial. Notice that (iii) is, in this case, equivalent to condition (\ref{cond allocation exists}). Therefore (i) $\Rightarrow$ (iii) follows from Theorem \ref{theo allocation exists}. What remains to be shown is that (iii) $\Rightarrow$ (ii).
To this aim, notice that when (iii) is verified and we consider the bipartite graph $\tilde{\mc G}=(\mc X\times \mc X, \tilde{\mc E})$ where $(x,y)\in\tilde{\mc E}$ iff $(x,y)\in E$, Hall's theorem guarantees the existence of a matching in $\tilde {\mc G}$ complete on the first set, namely a permutation $\sigma:\mc X\to\mc X$ such that $(x,\sigma(x))\in\tilde{\mc E}$ for every $x\in\mc X$. The corresponding state allocation $W^{\sigma}$ defined as in in (\ref{extremal Nash}) is a matching allocation state.
\end{proof}

We can now extend the result contained in Example \ref{ex: complete all}.

\begin{corollary}\label{cor: reg existence all} Suppose that $\mc G=(\mc X,\mc E)$ is any undirected regular graph and that $\alpha_x=a$, $\beta_x=b$ for every  $x\in\mc X$ with $a\leq b$. Then, there exists a (matching) allocation state.
\end{corollary}
\begin{proof} Let $s$ be the degree of each node in the graph. Fix any subset $D\subseteq \mc X$. If $\mc E_{D}$ is the set of directed edges starting from a node in $D$, we have that $s|D|=|\mc E_{D}|\geq s|N(D)|$. This implies that that $|D|\leq |N(D)|$. We conclude using Proposition \ref{prop: existence matching all}.
\end{proof}

Not necessarily a matching allocation state is the desirable one. In certain applications, security issues may rather require to fragment the data of each unit as much as possible. Suppose we are under the same assumptions than in previous result, namely $\mc G=(\mc X,\mc E)$ is an undirected regular graph with degree $s$, $\alpha_x=a$, $\beta_x=b$ for every  $x\in\mc X$ with $a\leq b$. If moreover $s$ divides $a$ we can also consider the 'diffused' allocation state given by
\beq\label{diffused allocation}W_{xy}=\frac{a}{s}A_{xy}\eeq
where $A$ is the adjacency matrix of $\mc G$. Notice that all these matrices $W$ can also be interpreted as valid allocation states for the case when the underlying graph is complete.

For graphs that are not regular, simple characterizations of the existence of allocations in general do not exist. However, sufficient conditions can be obtained as the result below shows and whose proof follows along the same line than the proof of Corollary \ref{cor: reg existence all}.

\begin{proposition} Let $\mc G$ be any graph with minimal out-degree $d_{min}$ and maximum in-degree $d^-_{max}$. Let
$a=\max_x\alpha_x$, $b=\min_x\beta_x$ and assume that $a\leq b d_{min}/d^-_{max}$. Then, there exists an allocation state.
\end{proposition}

The above result can not be improved: indeed in a star graph with $\alpha_x=a$, $\beta_x=b$ for every $x\in\mc X$, it is immediate to see that the condition $a\leq b d_{min}/d^-_{max}$ is  necessary for an allocation to exist.


%

\subsection{The optimal allocation problem}
On the set of allocation states, we define a reward functional measuring qualitative and realistic features that we desire the solution to possess, i.e., congestion and aggragation. Functionals considered in this paper have a separable structure that is a standard assumption in allocation problems \cite{DWG}.

We start with a notation. Given a (partial) allocation state $W\in\mc W_p$, we denote with the symbols $(W_{x\cdot})$ and $(W_{\cdot y})$, respectively, the row vector of W with label $x$, and the column vector of W with label $y$. We consider functionals $\Psi:\mc W_p\to\R$ of the type:
\begin{equation}\label{cost functional}\Psi(W)=\sum_{x\in\mc X}f_x(W_{x\cdot})+\sum_{y\in\mc X}g_y(W_{\cdot y})\end{equation}
consisting of two parts: one that takes into account the way each unit is succeeding in allocating its data and another that is typical a congestion term and considers the amount of data present in the various resources. Our goal is to maximize the functional $\Psi$ over the set of allocation states $\mc W$. The reason for defining $\Psi$ in the larger set of partial allocation states $\mc W_p$ will be clearer later when we present the game theoretic set up and the algorithm.

Examples and simulations in this paper will focus on the following cases:
\begin{equation}\label{cost functional example}\begin{array}{rcl}f_x(W_{x\cdot})&=&C^{all}\sum_y W_{xy}+C^{agg}\sum_{y\in\mc X}W_{xy}^2,\\[5pt] g_y(W_y)&=&-C^{con}_y(W_y)^2\end{array}\end{equation}
We now explain the meaning of the various terms:
\begin{itemize}
\item the term $C^{all}\sum_x\sum_y W_{xy}$ where $C^{all}> 0$ is sufficiently large, has the effect of pushing the optimum to be an allocation state (a configuration where all units have stored their entire set of data);
\item the term $C^{agg}\sum_x\sum_{y\in\mc X}W_{xy}^2$ has different significance depending on the sign of $C^{agg}$. If $C^{agg}>0$ plays the role of an aggregation term, it pushes units not to use many different resources for their allocation. If instead $C^{agg}<0$, the term has the opposite effect as it pushes towards fragmentation of the data.
\item the term $-\sum_yC^{con}_y(W_y)^2$ is a classical congestion term:  the constants $-C^{con}_y< 0$ for all $y$ measure the reliability of the various resources and pushes the use of more reliable resources
\end{itemize}
An alternative choice for the resource congestion term is the following. Put $|W_{\cdot y}|_H=|\{x\in \mc X,\, W_{xy}>0\}|$ the number of units that are using resource $y$ and consider
\beq\label{cost functional example 2} g_y(W_{\cdot y})=-C^{con}_y |W_{\cdot y}|_H\eeq
This might be useful in contexts where it is necessary the control the number of units accessing the same resource, to avoid communication burden. 

The functionals (\ref{cost functional example}) and (\ref{cost functional example 2}) reflects the features that we wanted to enforce: congestion and aggregation.
The reason for the latter feature comes from the fact that an exceeding fragmentation of the stored data will cause a blow up in the number of communications among the units, both in the storage and recovery phases. This feature should be considered against another feature, the diversification of back ups, which in this paper is not going to be addressed. In real applications, units will need to store multiple copies of their data in order to cope with security and failure phenomena. In that case, these multiple copies will need to be stored in different units. On the other hand, the congestion term is represented by a classical cost function that each user possibly experiences, for instance, as a delay in the storage/recovery actions.

The above desired features may be contradictory in general and we want to have tunable parameters to make the algorithm converge towards a desired compromised solution. The choice of this functionals has been made on the basis of simple realistic considerations and on the fact that, as exploited below, this leads to a potential game. In principle, different terms in the utility function can be introduced in order to make units to take into considerations other desired features (e.g multiple back up).

While our theory and algorithms will be formulated for a generic $\Psi$ as defined in (\ref{cost functional}), the example proposed and the numerical simulations will be restricted to the specific cases we have described.

Below we present a couple of examples of explicit computation of the maxima of $\Psi$. We assume $\Psi$ to be of the form described in (\ref{cost functional}) and (\ref{cost functional example}) with $C^{con}_y=C^{con}$ for every $y\in\mc X$. We also assume that $\alpha_x=a$, $\beta_x=b$ for every  $x\in\mc X$ with $a\leq b$.
\begin{example}\label{ex: homogeneous1}  Suppose that $\mc G=(\mc X,\mc E)$ is any undirected regular graph and assume that $\alpha_x=a$, $\beta_x=b$ for every  $x\in\mc X$ with $a\leq b$. Take $\Psi$ to be of the form described in (\ref{cost functional}) and (\ref{cost functional example}) with $C^{con}_y=C^{con}>0$ for every $y\in\mc X$. There are two cases:
\begin{itemize} 

\item $C^{agg}>0$. In this case, the maxima of $\Psi$ coincide with the matching allocation states. Indeed notice that any matching allocation state $W$ (whose existence is guaranteed by Proposition \ref{prop: existence matching all}) separately maximizes, for each $y$, the two expressions $\sum_{y\in\mc X}W_{xy}$ and $\sum_{y\in\mc X}W_{xy}^2$.  Moreover, considering that $W_y=a$ for every resource $y$, simultaneously, minimize the congestion expression $\sum_{y\in\mc X}(W_y)^2$. The fact that these are the only possible maxima is evident from these considerations.

\item  $C^{agg}<0$. If the degree $s$ of $\mc G$ divides $a$,  arguing like above, we see that the unique maximum is given by the diffused allocation state (\ref{diffused allocation}). When $s$ does not divide $a$, such a simple solution does not exist. In this case, maxima can be characterized as follows. Put $a=sk+r= (s-r)k+r(k+1)$ (with $r<s$) and consider a regular subgraph $\tilde G$ of degree $r$. An optimal allocation is obtained by letting units allocate $k+1$ atoms of their data in each of their neighbors in $\tilde G$ and $k$ atoms of their data in each of the remaining neighbors.
\end{itemize}

\end{example}

\section{The game theoretic set-up and the algorithm}\label{sec: algorithm}

In this paper we recast the optimization problem into a game theoretic context and we then use learning dynamics to derive decentralized algorithms adapted to the given graph topology that solve the allocation problem and maximize the functional $\Psi$.

Assume that a functional $\Psi$ as in (\ref{cost functional}) has been fixed. We associate a game to $\Psi$ according to the ideas developed in \cite{VTA,DWG} where there can be found other possible utility and potential functions. 

The set of actions $\mc A_x$ of a unit $x$ is given by all possible row vectors $(W_{x\cdot})$ such that $\sum_xW_{xy}\leq\alpha_x$. In this way the product set of actions $\prod_x\mc A_x$ can be made to coincide with the space of non-negative matrices $W\in\R^{\mc X\times \mc X}$ such that $\sum_xW_{xy}\leq\alpha_x$ for every $x\in\mc X$. Such a $W$ in general is not a partial allocation. Indeed, such a $W$ will automatically only possess properties (P1) and (P2). We have that $W\in\mc W_p$ if the extra conditions (P3), $\sum_yW_{xy}\leq\beta_y$ for every $y\in\mc X$, is satisfied. This is a key non classical feature of the game associated to our model: the storage limitations make the available actions of a unit depend on the choice made be the other ones.

Now, for each unit $x$, we define its utility function
$U_x:\mc W_p\to \RR$ as
\begin{equation}\label{utility1}U_x(W)=f_x(W_{x\cdot})+\sum_{y\in N_x^-}g_y(W_{\cdot y})\,.\end{equation}
Note that, in order to compute $U_x(W)$, unit $x$ needs to know, besides the state of its own data allocation $\{W_{x\cdot}\}$, the congestion state $g_y(W_{\cdot y})$ of the neighboring resources.

We now recall some basic facts of game theory. A \emph{Nash equilibrium}  is any allocation state $W\in\mc W_p$ such that, for every agent $\bar x\in\mc X$, and for every $W'\in\mc W_p$ such that $W_{xy}=W'_{xy}$ for every $x\neq \bar x$ and for every $y$, it holds
\begin{equation}\label{Nash} U_{\bar x}(W)\geq U_{\bar x}(W')
\end{equation} If $W, W'\in\mc W_p$ are two allocation states such that $W_{xy}=W'_{xy}$ for every $x\neq \bar x$ and for every $y$, it is straightforward to see that the following equality holds
\begin{equation}\label{potential1} U_x(W')-U_x(W)=\Psi(W')-\Psi(W)\end{equation}
This says, in the language of game theory \cite{Pot}, that the game is \emph{potential} with \emph{potential function} given by $\Psi$ itself. A simple classical result says that maxima of the potential are Nash equilibria for the game. In general the game will possess extra Nash equilibria. 

The choice (\ref{utility1}) is not the only one to lead to a potential game with potential $\Psi$. Other possibilities can be constructed following \cite{VTA,DWG}.
 As far as our theory is concerned, the specific form of the utility functions is not important as far as it leads to a potential game with potential $\Psi$. 
On the utility functions, (\ref{utility1}) or its possible alternatives, we impose a monotonicity condition that essentially says that no unit will ever have a vantage to remove data already allocated. Precisely, we assume that for every $W, W'\in \mc W_p$, $\bar x\in\mc X$ such that $W_{xy}=W'_{xy}$ for every $x\neq \bar x$ and for every $y$, the following holds
\beq\label{monotonicity} 
W'^{\bar x}<W^{\bar x}\;\Rightarrow\; U_{\bar x}(W')<U_{\bar x}(W)\eeq
This condition is not strictly necessary for our results (as our algorithm actually will not allow units to remove data), it is however a meaningful assumption and simulations show that helps to speed up the algorithm.

We now focus on the case when $\Psi$ is of the form given by (\ref{cost functional example}) with $C^{con}_y=C^{con}$ for every $y\in\mc X$. In this case, a simple check shows that the monotonicity condition
(\ref{monotonicity}) is guaranteed if we impose the condition
\beq\label{monotonicity2} C^{all}>2(||\alpha||_{\infty}|C^{agg}|+||\beta||_{\infty}C^{con})\eeq
where $||v||_{\infty}=\max v_i$ is the infinity norm of a vector.

We conclude this section, computing the Nash equilibria in a couple of simple examples and discussing the relation with the maxima of $\Psi$.
%
%
%
%
%
%
%
%

\begin{example} Suppose that $\mc G$ is the complete graph with three units and that $\alpha_x=a=2$ and $\beta_x=b\geq 2$ for $x=1,2,3$. Consider $\Psi$ to be of the form (\ref{cost functional example}) with $C^{con}_y=C^{con}$ for every $y\in\mc X$ and that condition (\ref{monotonicity2}) holds. Consider the following allocation states
$$W^1=\left(\begin{matrix}  0&2&0\\ 0&0&2\\ 2&0&0\end{matrix}\right),  W^2=\left(\begin{matrix}  0&0&2\\ 2&0&0\\ 0&2&0\end{matrix}\right), 
W^3=\left(\begin{matrix}  0&1&1\\ 1&0&1\\ 1&1&0\end{matrix}\right)$$
We know from the considerations in Example \ref{ex: homogeneous1} that in the case when $C^{agg}>0$, the matching allocation states $W^1$ and $W^2$ are the (only) two maxima of $\Psi$ and thus Nash equilibria. Instead, if $C^{agg}<0$, the diffused allocation state $W^3$ is the only maximum of $\Psi$ and is in this case a Nash equilibrium. 

Notice now that if $b<3$, the only three possible allocation states are $W^i$ for $i=1,2,3$. Since any two of these matrices differ in more than one row and condition (\ref{monotonicity2}) yields (\ref{monotonicity}), we deduce that all three of them are in this case Nash equilibria, independently on the sign of $C^{agg}$.

Suppose now that $b\geq 3$. Explicit simple computations show that, if $C^{agg}\leq C^{con}$, $W^3$ is a Nash equilibrium and if $C^{agg}\geq -6C^{con}$, $W^1$ and $W^2$ are Nash equilibria. In summary,
if $b<3$ or if $b\geq 3$ and $-6C^{con}\leq C^{agg}\leq C^{con}$, the three matrices $W^i$ for $i=1,2,3$ are Nash equilibria. 
\end{example}

%
%

The next example shows that also partial allocations may be Nash equilibria.

\begin{example}\label{partial allocation} $\mc G$ $5$-cycle, $\alpha_x=a=4$ and $\beta_1=7$, $\beta_2=2$, $\beta_3=4$, and $\beta_4=\beta_5=6$. It can be checked that the two matrices below are both Nash equilibria:
$$W=\left(\begin{matrix}  0&0&0&0&4\\ 3&0&1&0&0\\ 0&0&0&4&0\\ 0&0&3&0&1\\ 4&0&0&0&0\end{matrix}\right), \qquad W=\left(\begin{matrix}  0&0&0&0&4\\ 3&0&0&0&0\\ 0&0&0&4&0\\ 0&0&4&0&0\\ 4&0&0&0&0\end{matrix}\right)$$
The one of the right is a maximum of $\Psi$, the one on the left is instead a partial allocation.
\end{example}



\subsection{The algorithm}

The allocation algorithm we are proposing is fully distributed and asynchronous and is only based on communications between units, taking place along the links of the graph $\mc G=(\mc X, \mc E)$. It is based on the ideas of learning dynamics where, randomly, units activate and modify their action (allocation state) in order to increase their utility. The most popular of these dynamics is the so-called \emph{best response} where units at every step choose the action maximizing their utility. This dynamics is proven to converge almost surely, in finite time, to a Nash equilibrium. In presence of Nash equilibria that are not maxima of the potential (as it is in our case) best response dynamics is not guaranteed to converge to a maximum. This is simply because Nash Equilibria are always equilibrium points for the dynamics.
A popular variation of the best response is the so-called \emph{noisy} best response (also known as \emph{log-linear learning}) where maximization of utility is relaxed to a random choice dictated by a Gibbs probability distribution. 

We now illustrate the details of our algorithm. For the sake of proposing a realistic model we immagine that units may temporarily be shut down or in any case disconnected from the network. We model this assuming that, 
at every instant of time, a unit is either in functional state on or off: units in functional state off are not available for communication and for any action including storage and data retrieval. A unit, which is currently in state on, can activate and either newly allocate or move some data among the available resources (e.g. those neighbors that still have place available and that are on at that time). 
%
%
%
%
The functional state of the network at a certain time will be denoted by $\xi\in\{0,1\}^\mc X$: $\xi_x=1$ means that the unit $x$ is on. The times when units modify their functional state (off to on or on to off) and the times when units in functional state on activate are modeled as a family of independent Poisson clocks whose rates will be denoted (for unit $x$), respectively, $\nu_x^{on}$, $\nu_x^{off}$, and $\nu_x^{act}$. The functional state of the network as a function of time $\xi(t)$ is thus a continuous time Markov process whose components are independent Bernoulli processes. 

We now describe the core of the algorithm, namely the rules under which activated units can modify their allocation state.

We start with some notation.
Given a (possibly partial) allocation state $W\in\mc W_p$, a functional state $\xi\in\{0,1\}^{\mc X}$, and a unit $\bar x\in\mc X$ such that $\xi_{\bar x}=1$, define:
$$\mc W_{\bar x}(W,\xi)=\left\{W'\in\mc W_p:
\begin{aligned} 
&W'_{xy}=W_{xy} \text{ if } x\neq \bar x \text{ or } \xi_y=0\\
& W'^{\bar x}\geq W^{\bar x},\; W'\neq W
\end{aligned}\right\}.$$
$\mc W_{\bar x}(W,\xi)$ describes the possible partial allocation states obtainable from $W$ by modifications done by the unit $\bar x$: only the terms $W_{\bar xy}$ where $y$ is on can be modified and the total amount of allocated data $W^{\bar x}$ can only increase or remain equal. Since the sets $\mc W_{\bar x}(W,\xi)$ can in general be very large, it is convenient to consider the possibility that the algorithm might use a smaller set of actions where units either allocate new data or simply move data from one resource to another one. 

Given $(W,\xi)\in \mc W_p\times \{0,1\}^{\mc X}$ and a unit $\bar x$, define 
\beq N_{\bar x}(W,\xi):=\{y\in N_{\bar x}\,|\, W_{\bar x y}<\beta_y,\, \xi_y=1\}\eeq
the set of available neighbor resources for $\bar x$ under the allocation state $W$ and the functional state $\xi$: those that are on and still have space available.

A family of sets $\mc M_{\bar x}(W,\xi)\subseteq \mc W_{\bar x}(W,\xi)$, defined for each $\bar x\in\mc X$ and each $(W,\xi)\in \mc W_p\times \{0,1\}^{\mc X}$,  is called \textit{admissible} if
\begin{enumerate}
\item[(i)] $W^{\bar x}<\alpha_{\bar x},\, y\in N_{\bar x}(W,\xi),   \Rightarrow \exists n : W'=W+n e_{\bar xy}\in \mc M_{\bar x}(W,\xi)$;
\item[(ii)] $W_{\bar x y'}>0, \, \xi_{y'}=1,\, y''\in N_{\bar x}(W,\xi) \Rightarrow W'=W+ (e_{\bar xy''}-e_{\bar xy'})\in \mc M_{\bar x}(W,\xi)$;
\item[(iii)] $W'\in \mc M_{\bar x}(W,\xi)$ iff $W\in \mc M_{\bar x}(W',\xi)$ for every $W\in\mc W$.
\end{enumerate}

Conditions (i) and (ii) essentially asserts that when a unit has an available neighbor resource not yet saturated, then $\mc M_{\bar x}(W,\xi)$ must incorporate the possibility to newly allocate or transfer already allocated data into it. Condition (iii) instead simply says that when the functional state does not change and we are in an allocation state, any transformation can be reversed.
%

Examples of admissible families $\mc M_{\bar x}(W,\xi)$ are the following
\begin{enumerate}
\item $\mc M_{\bar x}(W,\xi)=\mc W_{\bar x}(W,\xi)$
\item $\mc M_{\bar x}(W,\xi)=\{W'\in\mc W_{\bar x}(W,\xi): \exists y, \exists n\; W'=W+ne_{\bar x y})\}
\cup \{W'\in\mc W_{\bar x}(W,\xi): \exists y',y'', \exists n\; W'=W+n(e_{\bar x y''}-e_{\bar x y'})\}$
\item $\mc M_{\bar x}(W,\xi)=\{W'\in\mc W_{\bar x}(W,\xi): \exists y, \exists n\in Q\; W'=W+ne_{\bar x y})\}
\cup \{W'\in\mc W_{\bar x}(W,\xi): \exists y',y'', \exists n\in Q\; W'=W+n(e_{\bar x y''}-e_{\bar x y'})\}$ where $Q\subseteq N$ and $1\in Q$.
\end{enumerate}
In the second case, modifications allowed are those where a unit either allocate a certain amount of new data into a single resource or it moves data from one resource to another one. The third case puts an extra constraint on the amount of data allocated or moved: the simplest case is $Q=\{1\}$, just an atomic piece of data is newly allocated or moved.
Simulation presented in this paper all fit in this third case with various possible sets $Q$.
%

Given an admissible family $\mc M_{\bar x}(W,\xi)$, we now define a \emph{Gibbs measure} on it as follows.
Given a parameter $\gamma>0$, put
$$Z^{(W,\xi)}_{\bar x}(\gamma)=\sum_{\tilde W\in \mc M_{\bar x}(W,\xi)}e^{\gamma U_{\bar x}(\tilde W)}$$
$$Z^{(W, W',\xi)}_{\bar x}(\gamma)=\max\left\{Z^{(W,\xi)}_{\bar x}(\gamma), Z^{(W',\xi)}_{\bar x}(\gamma)\right\}$$
Now define, for $W'\in \mc M_{\bar x}(W,\xi)$,
\begin{equation}\label{Gibbs}P^{(W,\xi)}_{\bar x}(W')=\left\{\begin{array}{ll}\frac{e^{\gamma U_{\bar x}(W')}}{Z^{(W,\xi)}_{\bar x}(\gamma)},\quad &{\rm if}\;  ||W||<||W'||\\[8pt]
\frac{e^{\gamma U_{\bar x}(W')}}{Z^{(W, W',\xi)}_{\bar x}(\gamma)},\quad &{\rm if}\;  ||W||=||W'||\end{array}\right.\end{equation}
where $||W||=\sum_{xy}W_{xy}$,
and complete it to a probability by putting
$$P^{(W,\xi)}_{\bar x}(W)=1-\sum_{W'\in \mc M_{\bar x}(W,\xi)}P^{(W,\xi)}_{\bar x}(W')$$


The algorithm is completely determined by the choice of the admissible family $\mc M_{\bar x}(W,\xi)$ and of the probabilities (\ref{Gibbs}). If unit $\bar x$ activates at time $t$, the systems is in partial allocation state $W(t)$, and in functional state $\xi(t)$, it will jump to the new partial allocation state $W'$ with probabilities given by
\begin{equation}\label{algorithm}P(W(t+)=W')=P^{(W(t),\xi(t))}_{\bar x}(W'),\; W'\in \mc M_{\bar x}(W(t),\xi(t))\end{equation}
If unit $\bar x$ chooses a $W'$ such that $||W'||>||W||$ we say that it makes an \emph{allocation move}, otherwise, if $||W'||=||W||$, we talk of a \emph{distribution move}.

\section{Analysis of the algorithm}\label{sec: analysis}
In this section we analyze the behavior of the algorithm introduced above. We will essentially show two results:
\begin{enumerate}
\item first, we prove that if the set of state allocation $\mc W$ is not empty (i.e. condition (\ref{cond allocation exists}) is satisfied), the algorithm above will reach such an allocation in bounded time with probability $1$ (e.g $W(t)\in\mc W$ for $t$ sufficiently large);
\item second, we show that, under a slightly stronger assumption than (\ref{cond allocation exists}), in the double limit $t\to +\infty$ and then $\gamma\to +\infty$, the process $W(t)$ induced by the algorithm will always converge, in law, to a Nash equilibrium that is a global maximum of the potential function $\Psi$. 
\end{enumerate}

In order to prove such results, it will be necessary to go through a number of intermediate technical steps. 

In the sequel we assume we have fixed a triple $(\mc G,\alpha, \beta)$ satisfying the existence condition (\ref{cond allocation exists}), an admissible family of sets $\mc M_{\bar x}(W,\xi)$ and we consider the allocation process $W(t)$ described by (\ref{algorithm}) with any possible initial condition $W(0)$.

By the way it has been defined, the process $W(t)$ is Markovian conditioned to the functional state process $\xi(t)$. If we consider the augmented process $(W(t), \xi(t))$, this is Markovian and its only non zero transition rates are described below:
\beq\label{transition1}\begin{array}{l}
\Lambda_{(W,\xi), (W,\xi')}=\left\{
\begin{array}{ll}
\nu_{\bar x}^{on}\quad &{\rm if}\, \xi_{\bar x}=0,\,\xi'_{\bar x}=1,\, \xi_x=\xi'_x\,\forall x\neq \bar x\\[2pt]
\nu_{\bar x}^{off}\quad &{\rm if}\, \xi_{\bar x}=1,\,\xi'_{\bar x}=0,\, \xi_x=\xi'_x\,\forall x\neq \bar x\end{array}
\right.\\[15pt]

%
\Lambda_{(W,\xi), (W', \xi)}=\nu_{\bar x}^{act}P^{(W,\xi)}_{\bar x}(W') \quad 
\hspace{0cm} {\rm if}\, \xi_{\bar x}=1,\, W'\in N_{\bar x}(W,\xi)
%
%
\end{array}
\eeq

We now introduce a graph on $\mc W_p$ that will be denoted by $\mc L_p$: an edge $(W,W')$ is present in $\mc L_p$ if and only if $W'\in\mc M_{\bar x}(W,\1)$.
Notice that, if $\nu_{\bar x}^{act}>0$ for every $\bar x$, this can be equivalently described as $\Lambda_{(W,\1),(W', \1)}>0$. 
The graph $\mc L_p$ thus describes the possible jumps of the process $W(t)$ conditioned to the fact that all resources are in functional state on. We want to stress the fact that the graph $\mc L_p$ depends on the triple $(\mc G,\alpha,\beta)$ as well on the choice of the admissible family $\mc M_{\bar x}(W,\1)$ but not on the particular choice of the functional $\Psi$ or of the utility functions $U_{\bar x}$.

Our strategy, in order to prove our first claim, will be to show that from any element $W\in\mc W_p$ there is a path in $\mc L_p$ to some element $W'\in \mc W$. 

Given $W\in \mc W_p$ we define the following subsets of units 
$$\mc X^f(W):=\{x\in\mc X\;|\; W^x=\alpha_x\},$$
$$\mc X^{sat}(W):=\{x\in\mc X\setminus\mc X^f(W)\;|\; \not\exists y\in N_x\;\hbox{\rm s.t}\; W_y< \beta_y\}$$
Units in $\mc X^f(W)$ are called {\emph{fully allocated}: these units have completed the allocation of their data under the state $W$. Units in $\mc X^{sat}(W)$ are called \emph{saturated}: they have not yet completed their allocation, however, under the current state $W$, they can not make any action, neither allocate, nor distribute. 
Finally, define
$$\mc W_p^{sat}:=\{W\in\mc W_p\setminus\mc W\;|\; \mc X=\mc X^f(W)\cup \mc X^{sat}(W)\}$$
It is clear that from any $W\in\mc W_p\setminus \mc W_p^{sat}$, there exists units that can make either an allocation or a distribution move. Instead, if we are in a state $W\in \mc W_p^{sat}$, there are units that are not fully allocated and all these units con not make any move. The only units that can possibly make a move are the fully allocated ones.
Notice that, because of condition (\ref{cond allocation exists}), for sure there exist resources $y$ such that $W_y<\beta_y$ and these resources are indeed exclusively connected to fully allocated units. The key point is to show that in a finite number of distribution moves, performed by fully allocated units, it is always possible to move some data atoms from resources connected to saturated units to resources with available space: this will then make possible a further allocation move.

For any fixed $W\in \mc W_p$, we can consider the following graph structure on $\mc X$ thought as set of resources: $\mc H_W=(\mc X,\mc E_W)$. Given $y_1,y_2\in\mc X$, there is an edge from $y_1$ to $y_2$  if and only if there exists $x\in\mc X$ for which
$$W_{xy_1}>0,\quad (x,y_2)\in\mc E
$$
The edge from $y_1$ to $y_2$ will be indicated with the symbol $y_1\to_x y_2$ (to also recall the unit $x$ involved). The presence of the edge means that the two resources $y_1$ and $y_2$ are in the neighborhood of a common unit $x$ that is using $y_1$ under $W$. This indicates that $x$ can in principle move some of its data currently stored in $y_1$ into resource $y_2$ if this last one is available. We have the following technical result

\begin{lemma}\label{lemma equilibrium 1} Suppose $(\mc G ,\alpha, \beta)$  satisfies (\ref{cond allocation exists}). Fix $W\in \mc W_p$ and let $\bar y\in\mc X$  be such that there exists $\bar x\in N_{\bar y}$ with $W^{\bar x}<\alpha_{\bar x}$.
Then, there exists a sequence
\beq\label{sequence} \bar y=y_0,\,x_0,\,y_1,\dots ,y_{t-1},\,x_{t-1},\,y_t\eeq
satisfying the following conditions
\begin{enumerate}
\item[(Sa)] Both families of the $y_k$'s  and of the $x_k$'s are each made of distinct elements;
\item[(Sb)] $y_k\to_{x_k} y_{k+1}$ for every $k=0,\dots ,t-1$;
\item[(Sc)] ${W}_{y_k}=\beta_{y_k}$ for every $k=0,\dots ,t-1$, and ${W}_{y_t}<\beta_{y_t}$.
\end{enumerate}
\end{lemma}
\begin{proof}
Let $\mc Y\subseteq \mc X$ be the subset of nodes that can be reached from $\bar y$ in $\mc H_W$.
Preliminarily, we prove that there exists $y'\in\mc Y$ such that $W_{y'}<\beta_{y'}$.
Let
$$\mc Z:=\{x\in\mc X\;|\; \exists y\in\mc Y,\, W_{xy}>0\}$$
and notice that, by the way $\mc Y$ and $\mc Z$ have been defined,
\beq\label{condZ} x\in\mc Z,\; (x,y)\in \mc E\;\Rightarrow\; y\in\mc Y\eeq
Suppose now that, contrarily to the thesis, $W_y\geq \beta_y$ for all $y\in\mc Y$. Then, 
\begin{equation}\label{estim}
\begin{aligned}
\sum\limits_{x\in \mc Z}\alpha_x &\leq \sum\limits_{y\in \mc Y}\beta_y\\
&=\sum\limits_{y\in \mc Y}W_y\\
&=\sum\limits_{y\in \mc Y}\sum\limits_{x\in\mc Z}W_{xy}\\
&=\sum\limits_{x\in\mc Z}W^x\\
&<\sum\limits_{x\in \mc Z}\alpha_x
\end{aligned}
\end{equation}
where the first inequality follows from (\ref{condZ}) and (\ref{cond allocation exists}), the first equality from the contradiction hypothesis, the second equality from the definition of $\mc Z$, the third equality again from (\ref{condZ}) and, finally, last inequality from the existence of $\bar x$. This is absurd and thus proves our claim.

Consider now a path of minimal length from $\bar y$ to $\mc Y$ in $\mc H_W$:
$$\bar y=y_0\to_{x_0}y_1\to_{x_1}\dots \to_{x_{t-2}} y_{t-1}\to_{x_{t-1}}y_t$$
and notice that the sequence $\bar y=y_0,\,x_0,\,y_1,\dots ,y_{t-1},\,x_{t-1},\,y_t$ will automatically satisfy properties (Sa) to (Sc).
\qed
\end{proof}

We are now ready to prove the first main result.

\begin{theorem} \label{theo main 1} Assume that the following conditions hold
\begin{enumerate}
\item $(\mc G ,\alpha, \beta)$  satisfies (\ref{cond allocation exists}). 
\item $\mc M_{\bar x}(W,\xi)$ is an admissible family.
\end{enumerate}
Then, for every $W\in\mc W_p$ there is a path in $\mc L_p$ to some element $W'\in \mc W$. 
\end{theorem}

\begin{proof}
We will prove the claim by a double induction process. To this aim we consider two indices associated to any $W\in\mc W_p\setminus \mc W$. The first one is defined by
$$m_W=\sum\limits_{x\in\mc X}(\alpha_x-W^x)\geq 1$$
To define the second, consider any  $\bar x\in \mc X\setminus \mc X^f(W)$. We can apply Lemma \ref{lemma equilibrium 1} to $W$ and any $\bar y\in N_{\bar x}$ and obtain that we can find a sequence of agents 
$\bar y=y_0,\,x_0,\,y_1,\dots ,\,y_{t-1},\,x_{t-1},\,y_t$
satisfying the properties (Sa), (Sb), and (Sc) above. Among all the possible choices of $\bar x\in\mc X$, $\bar y\in N_{\bar x}$ and of the corresponding sequence, assume we have chosen the one minimizing $t$ and denote such minimal $t$ by $t_W$. The induction process will be performed with respect to the lexicographic order induced by the pair $(m_W, t_W)$. 

In the case when $t_W=0$, it means we can find $\bar x\in\mc X$ and $\bar y\in\ N_{\bar x}$ such that $W_{\bar y}<\beta_{\bar y}$. This yields $\bar y\in N_{\bar x}(W,\1)$. Hence, by property (i) in the definition of an admissible family, it follows that there exists $n$ such that $W'=W+ne_{\bar x\bar y}\in \mc M_{\bar x}(W,\1)$. Notice that $m_{W'}<m_W$. In case $m_W=1$, this means that $W'\in\mc W$. 
%
%

Consider now any $W\in \mc W_p\setminus \mc W$ such that $t=t_W>1$. Let $\bar x\in\mc X$, $\bar y\in\ N_{\bar x}$  and the sequence $\bar y=y_0,\,x_0,\,y_1,\dots ,\,y_{t-1},\,x_{t-1},\,y_t$
satisfying the properties (Sa), (Sb), and (Sc) above. Since $W_{x_{t-1}y_{t-1}}>0$ and $y_t\in N_{x_{t-1}}(W,\1)$, it follows from property (ii) in the definition of admissible families that 
$W'=W-(e_{x_{t-1}y_{t-1}}-e_{x_{t-1}y_t})\in \mc M_{x_{t-1}}(W,\1)$. Since $W'_{y_{t-1}}<\beta_{y_{t-1}}$, for sure $t_{W'}<t_W$. The induction argument is thus complete. \qed

\end{proof}

\begin{corollary} \label{cor main 1} Consider the process $W(t)$ as defined in (\ref{algorithm}) and assume that the following conditions hold
\begin{enumerate}
\item $(\mc G ,\alpha, \beta)$  satisfies (\ref{cond allocation exists}). 
\item $\nu^{on}_x>0$ and $\nu^{act}_x>0$ for every $x\in\mc X$, 
\item $\mc M_{\bar x}(W,\xi)$ is an admissible family.
\end{enumerate}
Then, 
$$\P(\exists t_0\,|\, W(t)\in\mc W\,\forall t\geq t_0)=1$$
\end{corollary}
\begin{proof}
It follows from the form of the transition rates (\ref{transition1}) and assumption 2), that the process $(W(t), \xi(t))$, starting from any initial condition $(W, \xi)$, will reach $(W, \1)$ in bounded time with positive probability. Combining with Theorem \ref{theo main 1} and using again 2), it then follows that $(W(t), \xi(t))$ reaches a couple $(W', \1)$ for some $W'\in\mc W$ in bounded time with positive probability. Since, by definition of an admissible family, the set $\{(W',\xi),\, W'\in\mc W\}$ is invariant by the process $(W(t), \xi(t))$, standard results on Markov processes yield the thesis.
\end{proof}

We are now left with studying the process $W(t)$ on $\mc W$.  Noisy best response dynamics are known to yield reversible Markov processes. This is indeed the case also in our case once the process has reached the set of allocations $\mc W$. Precisely, the following result holds:

\begin{proposition}\label{prop time rev} Suppose that $\nu_x^{on}, \nu_x^{off}>0$ for all $x\in\mc X$. Then, $(W(t), \xi(t))$, restricted to $ \mc W\times \{0,1\}^\mc X$, is a time-reversible Markov process. More precisely, for every $(W,\xi),(W', \xi')\in\mc W\times \{0,1\}^\mc X$  it holds
\beq\label{time reversible}\rho_{(W,\xi)}\Lambda_{(W,\xi),(W',\xi')}=\rho_{(W',\xi')}\Lambda_{(W',\xi'),(W,\xi)}\eeq
where
\beq \rho_{(W,\xi)}=\left[\prod\limits_{x:\,\xi_x=1}\nu_x^{\rm on}\prod\limits_{x:\,\xi_x=0}\nu_x^{\rm off} \right]
e^{\gamma \Psi(W)}\eeq
\end{proposition}
\begin{proof}
It follows from relations (\ref{transition1}) and the definition of admissible families, that the only cases when $\Lambda_{(W,\xi),(W',\xi')}$ and $\Lambda_{(W',\xi'),(W,\xi)}$ are not both equal to zero are the following:
\begin{enumerate}
\item[(i)] $W'=W$, $\xi_{\bar x}=0,\,\xi'_{\bar x}=1,\, \xi_x=\xi'_x\,\forall x\neq \bar x$
\item[(ii)] $W'=W$, $\xi_{\bar x}=1,\,\xi'_{\bar x}=0,\, \xi_x=\xi'_x\,\forall x\neq \bar x$
\item[(iii)] $\xi'=\xi$, $W'\in N_{\bar x}(W,\xi)$.
\end{enumerate}
In case (i), we have that
$$\frac{ \rho_{(W,\xi)}}{ \rho_{(W,\xi')}}=\frac{\prod\limits_{x:\,\xi_x=1}\nu_x^{\rm on}\prod\limits_{x:\,\xi_x=0}\nu_x^{\rm off}}{\prod\limits_{x:\,\xi'_x=1}\nu_x^{\rm on}\prod\limits_{x:\,\xi'_x=0}\nu_x^{\rm off}}=\frac{\nu_{\bar x}^{\rm off}}{\nu_{\bar x}^{\rm on}}=
\frac{\Lambda_{(W,\xi'),(W,\xi)}}{\Lambda_{(W,\xi),(W,\xi')}}$$
Case (ii) can be analogously verified.
Consider now case (iii). 
Using relations (\ref{potential1}), (\ref{transition1}), and (\ref{Gibbs}), we obtain
$$\begin{array}{rcl}\frac{ \rho_{(W,\xi)}}{ \rho_{(W',\xi)}}&=&e^{\gamma (\Psi(W)-\Psi(W'))}=e^{\gamma ( U_{\bar x}(W)- U_{\bar x}(W'))}\\
&=&\frac{\Lambda_{(W',\xi),(W,\xi)}}{\Lambda_{(W,\xi),(W',\xi)}}\end{array}
$$
\qed
\end{proof}

We now show that under a slight stronger assumption than (\ref{cond allocation exists}), namely,
\beq\label{cond W ergodic}
\sum\limits_{x\in A}\alpha_x< \sum\limits_{y\in N(A)}\beta_y\quad\forall A\subseteq \mc X\,,
\eeq
the process $(W(t),\xi(t))$ restricted to $\mc W\times \{0,1\}^\mc X $ is ergodic. Denote by $\mc L$ the subgraph of $\mc L_p$ restricted to the set $\mc W$. Notice that, as a consequence of time-reversibility, $\mc L$ is an undirected graph. Ergodicity is equivalent to proving that $\mc L$ is connected. We start with a lemma analogous to previous Lemma \ref{lemma equilibrium 1}.

\begin{lemma}\label{lemma equilibrium 2} Suppose $(\mc G ,\alpha, \beta)$  satisfies (\ref{cond W ergodic}) and let $W\in \mc W$. Then, for every $\bar y\in\mc X$, there exists a sequence (\ref{sequence})
satisfying the conditions (Sa), (Sb), and (Sc) as in Lemma \ref{lemma equilibrium 1}. 
\end{lemma}
\begin{proof} It is sufficient to follow the steps of to the proof of Lemma \ref{lemma equilibrium 1} noticing that in (\ref{estim}) the first equality is now a strict inequality, while the last strict inequality becomes an equality. \qed
\end{proof}

If $W,W'\in\mc W$ are connected through a path in $\mc L$, we write that $W\sim W'$. Introduce the following distance on $\mc W$: if $W^1, W^2\in\mc W$
$$\delta(W^1,W^2)=\sum\limits_{x,y}|W^1_{xy}-W^2_{xy}|$$
A pair $\{W^1, W^2\}\in\mc W$ is said to be minimal if 
$$\delta (W^1,W^2)\leq \delta (W^{1'}, W^{2'})\;\;\forall W^{1'}\sim W^1,\; \forall W^{2'}\sim W^2$$
Notice that $\mc L$ is connected if and only if for any minimal pair $\{W^{1}, W^{2}\}$, it holds $W^{1}=W^{2}$. 

\begin{lemma}\label{lemma minimal1} Let $\{W^1, W^2\}$ be a minimal pair.
Suppose $y\in\mc X$ is such that $W^1_y<\beta_y$. Then, $W^1_{xy}=W^2_{xy}$ for all $x\in\mc X$.
\end{lemma}
\begin{proof}
Suppose by contradiction that $W^1_{xy}<W^2_{xy}$ for some $x\in\mc X$. Then, necessarily, there exists $y'\neq y$ such that $W^1_{xy'}>W^2_{xy'}$. Consider then $W^{1'}=W^1-e^{xy'}+e^{xy}$. Since $\delta (W^{1'},W^2)<\delta(W^1,W^2)$, this contradicts the minimality assumption. Thus $W^1_{xy}\geq W^2_{xy}$ for all $x\in\mc X$. This yields $W^2_y<\beta_y$. Exchanging the role of $W^1$ and $W^2$ we obtain the thesis. \qed
\end{proof}

\begin{proposition}\label{prop connected} If condition (\ref{cond W ergodic}) holds true, the graph $\mc L$ is connected.
\end{proposition}
\begin{proof}
Let $\{W^1, W^2\}$ be any minimal pair. We will prove that $W^1$ and $W^2$ are necessarily identical. Consider any resource $y$. It follows from Lemma \ref{lemma equilibrium 2} that we can find a sequence $y=y_0,\,x_0,\,y_1\cdots ,\,y_{t-1},\,x_{t-1},\,y_t$ satisfying the same (Sa), (Sb), and (Sc) with respect to the state allocation $W^1$. Among all the possible sequences, choose one with $t$ minimal for given $y$. We will prove by induction on $t$ that $W^1_{xy}=W^2_{xy}$ for all $x\in\mc X$.

If $t=0$, it means that $W^1_{y}<\beta_{y}$. It then follows from Lemma \ref{lemma minimal1} that $W^1_{xy}=W^2_{xy}$ for all $x\in\mc X$. Suppose now that the claim has been proven for all minimal pairs $\{W^1,W^2\}$ and any $y\in\mc X$ for which $t< \bar t$ (w.r. to $W^1$) and assume that $y=y_0,\,x_0,\,y_1\cdots ,\,y_{\bar t-1},\,x_{\bar t-1},\,y_{\bar t}$ satisfies the properties (Sa), (Sb), and (Sc) with respect to $W^1$. 

Since $W_{x_{\bar t-1}y_{\bar t-1}}>0$ and $y_{\bar t}\in N_{x_{\bar t-1}}(W^1,\1)$, it follows from property (ii) in the definition of admissible families that 
$W'^{1}=W^1-(e_{x_{\bar t-1}y_{\bar t-1}}-e_{x_{\bar t-1}y_{\bar t}})\in \mc M_{x_{\bar t-1}}(W^1,\1)$. In other words, $W'^{1}\sim W^1$.


Consider now $W^2$ and notice that Lemma \ref{lemma minimal1} yields $W^2_{x_{\bar t-1}y_{\bar t}}=W^1_{x_{\bar t-1}y_{\bar t}}<\beta_{y_{\bar t}}$. Define
$$W^{2'}=\left\{\begin{array}{ll}W^2\;&{\rm if}\, W^2_{x_{\bar t-1}y_{\bar t -1}}=0\\
W^2-e_{x_{\bar t-1}y_{\bar t-1}}+e_{x_{\bar t-1}y_{\bar t}}\;&{\rm if}\, W^2_{x_{\bar t-1}y_{\bar t -1}}>0\end{array}\right.$$
Again, by property (ii) in the definition of admissible families, it follows that $W'^{2}\sim W^2$.
Since $\delta(W^{1'},W^{2'})\leq \delta(W^{1},W^{2})$, this implies that also $\{W'^{1},W'^{2}\}$ is a minimal pair. Notice that $y=y_0,\,x_0,\,y_1\cdots ,\,y_{\bar t-1}$ satisfies (Sa), (Sb), and (Sc) with respect to $W^{1'}$. Therefore, by the induction hypotheses, it follows that $W^{1'}_{xy}=W^{2'}_{xy}$ for all $x\in\mc X$. Since $W^1_{xy}=W^{1'}_{xy}$ and $W^2_{xy}=W^{2'}_{xy}$, result follows immediately.
\qed
\end{proof}

We can now state our final result.
\begin{corollary}\label{cor main 2} Assume that the following conditions hold
\begin{enumerate}
\item $(\mc G ,\alpha, \beta)$  satisfies (\ref{cond W ergodic}). 
\item $\nu^{on}_x>0$,  $\nu^{off}_x>0$, and $\nu^{act}_x>0$ for every $x\in\mc X$, 
\item $\mc M_{\bar x}(W,\xi)$ is an admissible family.
\end{enumerate}
Then, $(W(t),\xi(t))$, restricted to $\mc W\times \{0,1\}^{\mc X}$, is an ergodic time-reversible Markov process whose unique invariant probability measure is given by
$$\mu_{\gamma}(\xi, W)=Z_\gamma^{-1}\left[\prod\limits_{x:\,\xi_x=1}\nu_x^{\rm on}\prod\limits_{x:\,\xi_x=0}\nu_x^{\rm off} \right]
e^{\gamma \Psi(W)}$$
where $Z_{\gamma}$ is the normalizing constant.
\end{corollary}

\begin{proof}
Let $(W, \xi), (W', \xi')\in \mc W\times \{0,1\}^{\mc X}$. 
It follows from the form of the transition rates (\ref{transition1}) and the fact that $\nu^{on}_x>0$ for all $x$, that the process $(W(t), \xi(t))$, starting from $(W, \xi)$, will reach $(W, \1)$ in bounded time with positive probability. Combining with Proposition  \ref{prop connected} and using the fact that $\nu^{act}_x>0$ for all $x$, it then follows that $(W(t), \xi(t))$ reaches $(W', \1)$ in bounded time with positive probability. Finally, from $(W', \1)$ again the process reaches $(W', \xi')$ in bounded time with positive probability. This says that the process is ergodic and it thus possesses a unique invariant measure whose form can be derived by the time-reversibility property characterized in Proposition \ref{prop time rev}.\qed
\end{proof}

\medskip
{\bf Remark:} It follows from previous result that the process $W(t)$ converges in law to the probability distribution
$$\tilde\mu_{\gamma}(W):=\tilde Z_\gamma^{-1}e^{\gamma \Psi(W)}$$
Notice that when $\gamma\to +\infty$, the probability $\tilde\mu_{\gamma}$ converges to a probability concentrated on the set $\argmax_{W\in\mc W} \Psi(W)$ of state allocations maximizing the potential. Thus, if $\gamma$ is large, the distribution of the process $W(t)$ for $t$ sufficiently large will be close to a maximum of $\Psi$.

\medskip
{\bf Remark:} Condition (\ref{cond W ergodic}) is necessary for ergodicity. Notice indeed that in the case when $\mc G$ is complete and $\alpha_x=\beta_x=a$ for all $x\in\mc X$, under every allocation $W$ such that $W_y=a$ for every $y$, all resources will be saturated and, consequently, no distribution move will be allowed in $W$. Such allocations $W$ are thus all sinks in the graph $\mc L$ that is therefore not connected.

\section{Simulation}\label{sec simulation}
In this section, we present some numerical results to validate the theoretical approach. We aim to show that the algorithm is feasible for practical implementation and that it has good performance and scaling properties. Example presented are admittedly simple: our goal is here is not to work out codes with optimized performance, neither to present exhaustive sets of simulations.

All our examples are for the  case when the functional has the form defined in (\ref{cost functional example}), but the last one where instead we consider the form (\ref{cost functional example 2}).

We always take 
\begin{equation}\label{c-all}
C^{all}= 3(||\alpha||_{\infty}|C^{agg}| +||\beta||_{\infty}C^{con})
\end{equation}
This choice is motivated by the considerations in (\ref{monotonicity}). 
We assume that $C^{con}_y=1$ for all units and we consider both the case when the aggregation parameter $C^{agg}$ is positive or negative. 

As a graph, we use either a complete graph or a regular graph of degree $10$ randomly constructed according to the classical configuration model. 

We assume the admissible family $\mc M_{\bar x}(W,\xi)$ to be of type 3) presented before where
modifications allowed are those where a unit either allocate or move a number of data constrained in a set $Q$. Most of the examples are for $Q=\{1\}$: just one data is allocated or moved each time.

On the basis of our theoretical analysis, the algorithm, in the limit when $t\to +\infty$ and the inverse noisy parameter $\gamma \to +\infty$, is known to converge to the optimum. In practical implementations, a typical choice in these cases is to take the parameter $\gamma$, time-varying and diverging to $+\infty$. The tuning of the divergence rate is known to be critical to obtain good results. Here we have chosen the activation rate $\nu_{act}=1/n$ and 
$$\gamma(t+1)=\gamma(t)+\frac{1}{100000}$$
Moreover, we suppose the units to be always on (otherwise things get simply slowed down). The time horizon is fixed $T=5*\sum_{x\in\mc X}\alpha_x$: in this way a unit $x$ will activate, on average, a number of times equal to 5 times the number of data  it needs to allocate. As we will see, this time range is sufficient for the allocation to be completed and to get close to the optimum (this has been checked in those cases when the optimum is analytically known).

%
%

For all examples, the performance of the algorithm is analyzed considering the following parameters computed, in a Montecarlo style, by averaging over 10 runs of the algorithm. 
\begin{itemize}
 \item {\bf{Distance from Full Allocation:}} 
 $$\Delta=\sum_{x\in\mc X}\alpha_x-\sum_{x,y\in\mc X}W_{xy}$$
counts the quantity of atoms not yet allocated.
If the allocation is complete, this parameter is $0$.
\item {\bf Allocation complexity:} Denoted by $m_x$  the total number of moves (allocation and distribution) made by unit $x$, we consider
$$\nu_{moves}=\frac{1}{n}\sum\limits_{x\in\mc X}\frac{m_x}{\alpha_x}$$
$\nu_{moves}$ measures the number of allocation and distribution moves per piece of data. Since moving data from a resource to another can be expensive, it is an interesting parameter to consider

\item {\bf Distance in ratio from optimum:} In cases when the maximum of the potential $\Psi_{opt}$ is explicitly known (Example \ref{ex: homogeneous1}), we consider $\psi=\frac{\Psi(W^{T)}}{\Psi_{opt}}$. 
%
\item{\bf Degree  complexity:} We consider the average number of resources used by a unit
$$d:=\frac{1}{n}\sum_{x\in\mc X}\sum_{y\in\mc X}\1_{\{W^T_{xy}>0\}}$$ 
$d$ is a measure of how concentrated or diffused is the allocation. For matching allocations $d=1$. 
\end{itemize}

We now present a number of simulations for the case when the functional has the form defined in (\ref{cost functional example}).
We first consider the case when $Q=\{1\}$: just one data is allocated or moved each time a unit activates. We always take $C^{con}=1$ and $C^{all}$ chosen according to (\ref{c-all}) and different values for $C^{agg}$.

\begin{example}\label{10user}
Consider to have $n=10$ users on a complete graph such that $\alpha_x=a=45$ and $\beta_x=b=50$ for every unit $x$. 
We consider the cases: $C^{agg}=-7,-1,1/2,3$.
First, we show the final states reached by the dynamics for a single run of the algorithm  
$$W_{-7}=\left(\begin{array}{cccccccccc} 
0 & 5 & 5 & 5 & 5 & 5 & 5 & 5 & 5 & 5\\ 
5 & 0 & 5 & 5 & 5 & 5 & 5 & 5 & 5 & 5\\ 
5 & 5 & 0 & 5 & 5 & 5 & 5 & 5 & 5 & 5\\ 
6 & 5 & 5 & 0 & 6 & 5 & 6 & 6 & 0 & 6\\ 
5 & 5 & 5 & 5 & 0 & 5 & 5 & 5 & 5 & 5\\ 
5 & 5 & 5 & 5 & 5 & 0 & 5 & 5 & 5 & 5\\ 
5 & 5 & 5 & 5 & 5 & 5 & 0 & 5 & 5 & 5\\ 
5 & 5 & 5 & 5 & 5 & 5 & 5 & 0 & 5 & 5\\ 
5 & 5 & 5 & 5 & 5 & 5 & 5 & 5 & 0 & 5\\ 
5 & 5 & 5 & 5 & 5 & 5 & 5 & 5 & 5 & 0 
\end{array}\right)$$

$$W_{-1}=\left(\begin{array}{cccccccccc} 
     0&     5&     3 &    5   &  5   &  4  &   8&     7   &  6   &  2\\
     9  &   0  &   3   &  4     &4   &  1     &6&     7     &6 &    5\\
     1    & 3    & 0     &6   &  5     &6   &  6  &   3   &  8   &  7\\
     6&     1&     6&     0&     6&    10&     7   &  0&     5    & 4\\
     3  &   6  &   7  &   8  &   0    & 6    & 2    & 8&     0    & 5\\
     2    & 8    & 9    & 6    & 4  &   0   &  3    & 7  &   5     &1\\
     8     &3&     7&     5&     9&     7  &   0  &   0    & 0     &6\\
     4    & 5  &   4  &   4     &3  &   4   &  5  &   0   &  8&     8\\
     1  &   9    & 2    & 3   &  8 &    1    & 7   &  8  &   0  &   6\\
     8&    10     &5     &3&     3  &   7 &    3   &  4   &  2   &  0\\
\end{array}\right)$$

$$W_{1/2}=\left(\begin{array}{cccccccccc} 
     0&     0&    23&     7&     0&     4&     0&     0 &   10&     1\\
     3  &   0  &  12  &   7  &   2  &   7  &   2  &   5   &  0   &  7\\
     1    & 4    & 0     &2    &15    & 2    & 7    & 1    &12    & 1\\
    15     &6     &2     &0     &2    &14     &0     &1     &4     &1\\
     4&     4 &    0&     1&     0&     4&    13&     0&     3&    16\\
    18 &    0  &   0 &    1 &    3 &    0 &    8  &   2  &   8  &   5\\
     2    & 7     &6    &12   &  2   &  3   &  0    & 0    & 5    & 8\\
     1&     4&     1&     9&    12&     1&     9&     0&     4&     4\\
     0  &  17 &    1 &    1 &    4  &   1  &   0  &  16 &    0 &    5\\
     2    & 5    & 0    & 1    & 6    &10    & 2    &17   &  2   &  0\\
\end{array}\right)$$

$$W_{3}=\left(\begin{array}{cccccccccc} 
0 & 0 & 45 & 0 & 0 & 0 & 0 & 0 & 0 & 0\\ 
45 & 0 & 0 & 0 & 0 & 0 & 0 & 0 & 0 & 0\\ 
0 & 0 & 0 & 0 & 0 & 0 & 0 & 0 & 45 & 0\\ 
0 & 0 & 0 & 0 & 0 & 0 & 45 & 0 & 0 & 0\\ 
0 & 0 & 0 & 45 & 0 & 0 & 0 & 0 & 0 & 0\\ 
0 & 45 & 0 & 0 & 0 & 0 & 0 & 0 & 0 & 0\\ 
0 & 0 & 0 & 0 & 0 & 0 & 0 & 45 & 0 & 0\\ 
0 & 0 & 0 & 0 & 0 & 45 & 0 & 0 & 0 & 0\\ 
0 & 0 & 0 & 0 & 0 & 0 & 0 & 0 & 0 & 45\\ 
0 & 0 & 0 & 0 & 45 & 0 & 0 & 0 & 0 & 0 
\end{array}\right)$$

For the same runs, we plot in Figure 1  the time evolution of the potentials and we confront it with the optimal potential represented by the red line. For $C^{agg}=3$ a matching allocation state is reached and it is a maximum of $\Psi$ in this case. For $C^{agg}=-7$ the solution is also very close to the maximum that is the diffused allocation state. For $C^{agg}=1/2,-1$, the presence of Nash equilibria that are not maxima of $\Psi$ slows down the dynamics and the algorithm does not reach the maximum at time $T$ (this particularly evident for the case $C^{agg}=1/2$). 
%
Increasing in this case the time horizon to $T=20*\sum\alpha$, the final state of the system gets quite close to the maximum as confirmed by the two plots in Figure 2.

%
%
%

\begin{center}
\begin{figure}\label{fig_pot_t5}
\centering %
\subfigure[$C^{agg}=-7$]
{%
\includegraphics[scale=.2]{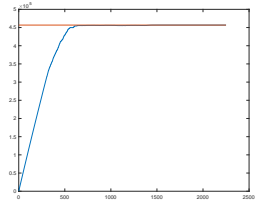}
}%
\hspace{1cm}
\subfigure[$C^{agg}=-1$]
{%
\includegraphics[scale=.2]{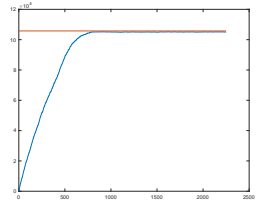}
}

\subfigure[$C^{agg}=1/2$]
{%
\includegraphics[scale=.2]{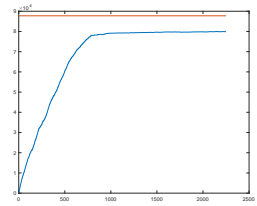}
}%
\hspace{1cm}
\subfigure[$C^{agg}=3$]
{%
\includegraphics[scale=.2]{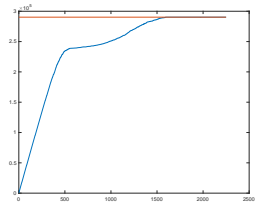}
}
\caption{Time evolution of the Potential
}
\end{figure}
\end{center}

\begin{center}
\begin{figure}[h]\label{fig_pot_t20}
\centering %
\subfigure[$C^{agg}=-1$]
{%
\includegraphics[scale=.2]{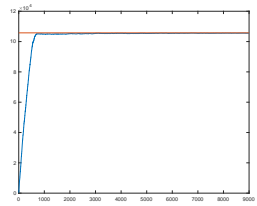}
}%
\hspace{1cm}
\subfigure[$C^{agg}=1/2$]
{%
\includegraphics[scale=.2]{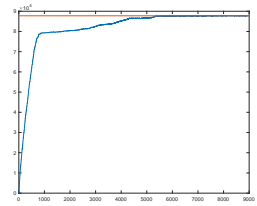}
}
\caption{Evolution of the Potential with $T=20*\sum_{x\in\mc X}\alpha_x$}
\end{figure}
\end{center}

The following table shows the performance parameters in the Montecarlo simulation for the usual $T$.
\begin{center}
\begin{table}[htb]\label{tab_esempio1}
\caption{Performance parameters for $n=10$}
\centering
\begin{tabular}{|l|c|c|c|c|}
\hline
& $C^{agg}=-7$ & $C^{agg}=-1$ & $C^{agg}=1/2$ & $C^{agg}=3$\\
\hline
$d$  & 9   & 8.7400  & 6.6200 & 1  \\[1pt]
$\psi$  &  1  &   0.9944   & 0.9156 & 1  \\
$\Delta$  & 0 & 0 & 0  & 0 \\
$\nu_{moves}$ & 3.1669 &  4.9389 &  4.9331   & 3.2449 \\
\hline
\end{tabular}
\end{table}
\end{center}   
%
\end{example}

\newpage
From now on we focus on the cases $C^{agg}=-7,3$, $C^{con}=1$ and $C^{all}$ chosen according to (\ref{c-all}), showing that reasonably good properties are maintained for larger communities and different topologies. 

\begin{example}
Consider to have $n=50$ users on a complete graph and on a regular graph of degree $10$ such that $\alpha_x=a=45$ and $\beta_x=b=50$ for every unit $x$. 
While a matching allocation state is not reached when $C^{agg}=3$, the value of the average degree shows that the solution is quite concentrated with most of units allocating in just one resource. Instead, for $C^{agg}=-7$ we have reached an optimum diffused allocation.
\begin{center}
\begin{table}[h]
\caption{Performance parameters for $n=50$}
\centering
\subtable[$C^{agg}=3$]{
\begin{tabular}{|l|c|c|}
\hline
& Complete & Regular  \\
\hline
$d$ & 1.2400  & 1.2280  \\[1pt]
$\psi$ &   0.9794 & 0.9872 \\
$\Delta$ & 0 & 0\\
$\nu_{moves}$ &  1.8238 & 2.4538   \\
\hline
\end{tabular}}
%
\subtable[$C^{agg}=-7$]{
\begin{tabular}{|l|c|c|}
\hline
& Complete & Regular  \\
\hline
$d$ &  45 & 10  \\[1pt]
$\psi$ &  1  & 1  \\
$\Delta$ & 0 & 0\\
$\nu_{moves}$ &  1.3746 & 1.2898   \\
\hline
\end{tabular}}
\end{table}
\end{center} 

\end{example}
The next example shows how the presence of heterogeneous resources does not alter much the performance of the algorithm

\begin{example}
Consider to have $n=50$ users on a complete graph and on a regular graph of degree $10$ such that $\alpha_x=a=43$ for every $x$. Assume that half of the units have $\beta_x=40$ and half or them instead $\beta_x=50$. Notice that, in this case, for the regular graph topology, there is no a-priori guarantee that allocation is feasible. Simulations show however that allocation is reached in all cases. 

 
%
\begin{center}
\begin{table}[h]
\caption{Performance parameters for $n=50$, resources with different storage capabilities}
\centering
\subtable[$C^{agg}=3$]{
\begin{tabular}{|l|c|c|}
\hline
& Complete & Regular\\
\hline
$d$ &  2.0040 & 2.2760 \\[1pt]
$\Delta$ & 0 & 0\\
$\nu_{moves}$ & 2.1540  & 4.1273 \\
\hline
\end{tabular}}
%
\subtable[$C^{agg}=-7$]{
\begin{tabular}{|l|c|c|}
\hline
& Complete & Regular\\
\hline
$d$ & 42.628  & 10 \\[1pt]
$\Delta$ & 0 & 0\\
$\nu_{moves}$ &  1.9754  &  1.2862\\
\hline
\end{tabular}}
\end{table}
\end{center} 

\end{example}

In the following example we consider larger families of units connected through a regular graph of degree $10$. Numerical results show the good scalability properties of the algorithm.
\begin{example}
Suppose to have $n=100,200,300$ users on a regular graph of degree 10 with $\alpha_x=a=45$ and $\beta_x=b=50$. Table 4 shows the performance parameters. 
\begin{center}
\begin{table}[h]
\caption{Performance parameters for $n=100,100,300$}
\centering
\subfigure[$C^{agg}=-7$]{
\centering
\begin{tabular}{|l|c|c|c|}
\hline
& $n=100$ & $n=200$ & $n=300$\\
\hline
$d$ &  10 &  10 & 10 \\[1pt]
$\psi$ &   1 & 1 & 1 \\
$\Delta$ & 0 & 0 & 0 \\
$\nu_{moves}$ &  1.2535  & 1.2832  &  1.2897  \\
\hline
\end{tabular}}
\centering
\subfigure[$C^{agg}=3$]{
\begin{tabular}{|l|c|c|c|}
\hline
& $n=100$ & $n=200$ & $n=300$\\
\hline
$d$ &  1.3980 &  1.4040 & 1.4017 \\[1pt]
$\psi$ &   0.9751 & 0.9753 & 0.9748 \\
$\Delta$ & 0 & 0 & 0 \\
$\nu_{moves}$ &  2.0125  & 1.6346  &  1.5114  \\
\hline
\end{tabular}}

\end{table}
\end{center} 

\end{example}

Next example consider the case when allocations and distributions are allowed with different granularity $Q$. 

\begin{example}
Consider to have $n=10$ users on a complete graph such that $\alpha_x=a=45$ and $\beta_x=b=50$ for every unit $x$. 
We assume that units can allocate or move each time a quantity of data belonging to either $Q_1=\{1,5,10\}$ or $Q_2=\{1,25,45\}$. We also report the case $Q_0=\{1\}$ for the sake of comparison.
%
\begin{center}
\begin{table}[h]
\caption{Performance parameters for $n=10$, different granularity}
\centering
\subtable[$C^{agg}=-7$]{
\begin{tabular}{|l|c|c|c|}
\hline
& $Q_0$ & $Q_1$ & $Q_2$   \\
\hline
$d$ & 9  &  9 &  3.8500    \\[1pt]
$\psi$   & 1 &   0.9999 & 0.8902     \\
$\Delta$  & 0  & 0 & 0    \\
$\nu_{moves}$  & 3.1669 &  0.2311  &     0.1767       \\
\hline
\end{tabular}}
%
\subtable[$C^{agg}=3$]{
\begin{tabular}{|l|c|c|c|}
\hline
& $Q_0$ & $Q_1$ & $Q_2$ \\
\hline
$d$   & 1 &  1.0100 &  1.0100   \\[1pt]
$\psi$   & 1 &   0.9996 & 0.9999   \\
$\Delta$  & 0 & 0 & 0    \\
$\nu_{moves}$  & 3.2449 &  0.1224  &     0.0229    \\
\hline
\end{tabular}}
\end{table}
\end{center} 
As expected, the possibility to allocate at one time larger sets of data drastically reduces the number of allocation and distribution moves and speeds up the algorithm. Notice however that in one case, using the set $Q_2$, the algorithm does not reach the maximum. This phenomenon is probably due to the fact that  allocating large set of data at once can lead to allocation states quite far from the optimum and thus require longer time to converge. This says that the choice of the set $Q$ is likely to play a crucial role in order to optimize the speed of convergence of the algorithm.

\end{example}

Finally, the last example is for the objective functional with the alternative congestion term (\ref{cost functional example 2}).


\begin{example}
Consider to have $n=50$ users on a complete graph and on regular graph of degree $10$ such that $\alpha_x=a=45$ and $\beta_x=b=50$ for every unit $x$. We take $C^{con}=1$ and $C^{all}$ chosen according to (\ref{c-all}) while we take different values for $C^{agg}$.
As expected, in this case, varying the aggregation parameter $C^{agg}$, we obtain solutions with a different degree of fragmentation. 
\begin{center}
\begin{table}[h]
\caption{Performance parameters for $n=50$, alternative congestion term}
\centering
\subtable[Regular graph]{
\begin{tabular}{|l|c|c|c|c|c|c|c|}
\hline
&  $C^{agg}=0$&$C^{agg}=-1/2$&$C^{agg}=-2$&   $C^{agg}=-5$    \\
\hline
$d$ & 4.7680 & 8.2849  &  9.8240 &    10    \\[1pt]
$\Delta$ & 0 &0 & 0 &0 \\
$\nu_{moves}$ & 4.8604 & 4.4729  & 4.3748    & 5   \\
\hline
\end{tabular}}
%
%
\subtable[Complete graph]{
\begin{tabular}{|l|c|c|c|c|c|c|c|}
\hline
&   $C^{agg}=-2$  & $C^{agg}=-10$ & $C^{agg}=-20$ & $C^{agg}=-100$  \\
\hline
$d$  &    13.9640    & 21.6400 & 23 & 45   \\[1pt]
$\Delta$   & 0&0 &0 &0  \\
$\nu_{moves}$  &  4.7415    & 4.9951 & 5 & 1.3547  \\
\hline
\end{tabular}}
\end{table}
\end{center}  
This is particularly evident in the case of a complete graph. The choice of the topology and of the functional parameters can be seen, in this case, as alternative or complementary ways to prescribe the complexity of the allocation in terms of links used.

\end{example}

\section{Conclusions}
We have presented and mathematically analyzed a decentralized allocation algorithm, motivated by the recent interest in cooperative cloud storage models. We have proved convergence and we have shown the practical implementability of the algorithm. The tuning of its parameters to optimize performance will be considered elsewhere. In this direction, it will also be useful to investigate the possibility to use different utility functions in the definition of the algorithm, following the ideas in \cite{DGDO} and \cite{DWG}. On the other hand, it would be of interest to deepen the relation of this model with generalized Nash equilibrium problems \cite{GNEP,GenPot} to better understand the level of generality of out approach.

\section*{Acknowledgment}
We acknowledge that this work has been done while Barbara Franci was a PhD Student sponsored by a TIM/Telecom Italia grant.

\end{document}